\documentclass[12pt]{amsart}

\usepackage[all]{xy}
\usepackage{mathrsfs}
\usepackage{amsmath}
\usepackage{latexsym}
\usepackage{amssymb}
\usepackage{amsfonts}
\usepackage{longtable}
\usepackage{colortbl}
\include{PDF}
\usepackage{graphics}

\def\l{\langle} \def\r{\rangle} 
 
 \def\ZZ{\mathbb Z}

\def\Aut{{\sf Aut}} 

 \def\K{{\sf K}}

 \def\soc{{\sf soc}}

\def\D{{\rm D}} 
\def\S{{\rm S}} 
 \def\M{{\rm M}}

\def\C{{\bf C}}\def\N{{\bf N}}
\def\O{{\bf O}}

\def\calO{{\mathcal O}}

\def\Ome{{\it \Omega}}
\def\Ga{{\it \Gamma}}

\def\Del{{\it\Delta}}
\def\a{\alpha} \def\b{\beta} \def\g{\gamma} \def\s{\sigma}
\def\t{\tau}

 \def\GL{{\rm GL}}

\def\Sp{{\rm Sp}}

\def\GammaL{{\rm \Gamma L}}

\def\PGammaL{{\rm P\Gamma L}}

\def\SigmaL{{\rm \Sigma L}}

\def\A{{\rm A}}
\def\Sym{{\rm Sym}}

\def\PSL{{\rm PSL}}  \def\PGL{{\rm PGL}}
\def\GL{{\rm GL}} \def\SL{{\rm SL}}
\def\AGL{{\rm AGL}}

  \def\D{{\rm D}}

\def\Alt{{\rm Alt}}

\def\K{{\bf K}}

\newtheorem{thm}{Theorem}[section]

\newtheorem{lem}[thm]{Lemma}

\newtheorem*{remark}{Remark}

\newtheorem{lemma}[thm]{Lemma}
\newtheorem{corollary}[thm]{Corollary}
\newtheorem{example}[thm]{Example}
\newtheorem{theorem}{Theorem}[section]

\def\qed{{\hfill$\Box$\smallskip}
\medbreak}

\begin{document}

\title[Alternating groups and 2-arc-transitive graphs]
%{Subgroups of alternating groups of odd indices, and $2$-arc-transitive graphs}
{Two-arc-transitive graphs of odd order -- II}
\thanks{The project was partially supported by the NNSF of China
(11771200, 11931005, 11861012,  11971248, 11731002) and the Fundamental Research Funds for the Central Universities.}

\author[Li]{Cai Heng Li}
\address{Cai Heng Li\\
Department of Mathematics \\
Southern University of Science and Technoogy\\
Shenzhen 518055, P. R. China}
\email{lich@sustech.edu.cn}

\author[Li]{Jing Jian Li}
\address{Jing Jian Li\\ School of Mathematics and Information Sciences\\ Guangxi
University\\ Nanning 530004, P. R. China.}
\email{lijjhx@gxu.edu.cn}

\author[Lu]{Zai Ping Lu}
\address{Zai Ping Lu\\ Center for Combinatorics,
 LPMC\\ Nankai University\\
 Tianjin 300071,
P. R. China} \email{lu@nankai.edu.cn}

\begin{abstract}
It is shown that each subgroup of odd index in an alternating group of degree at least 10 has all insoluble composition factors to be alternating.
A classification is then given of 2-arc-transitive graphs of odd order admitting an alternating group or a symmetric group.
This is the second   of a series of papers aiming towards a classification of 2-arc-transitive graphs of odd order.
\end{abstract}

\maketitle

\date\today

\vskip 20pt

\section{ Introduction}

%The principal motivation of discovering the result was from algebraic graph theory.
Let $\Ga=(V,E)$ be a graph with vertex set $V$ and edge set $E$, which is finite, simple and undirected.
The number of vertices $|V|$ is called the {\it order} of the graph.
A $2$-{\it arc} in $\Ga$ is a triple of distinct vertices $(\a,\b,\g)$
such that $\b$ is adjacent to both $\a$ and $\g$.
In general, for an integer $s\geqslant  1$, an $s$-{\it arc} is a sequence of $s+1$ vertices with any two consecutive vertices adjacent and any three consecutive vertices distinct.
A graph $\Ga$ is said to be {\it $(G,s)$-arc-transitive} if $G\leqslant\Aut\Ga$ is transitive on both the vertex set and the set $s$-arcs of $\Ga$, or simply called {\it $s$-arc-transitive}.
By the definition,
an $s$-arc-transitive graph is also $t$-arc-transitive for $1\leqslant  t<s$.

The class of $s$-arc-transitive graphs has been one of the central topics in algebraic graph theory since Tutte's seminal result \cite{Tutte}: there is no 6-arc-transitive cubic graph,
refer to \cite{Trofimov, Weiss} and \cite{Baddeley,FP2,FP1,HNP,IP,Li3,LSS,Li4,Prag-o'Nan}, and references therein.
A great achievement in the area was due to Weiss \cite{Weiss} who proved that there is no 8-arc-transitive graph of valency at least 3.
Later in \cite{Li1}, the first named author proved that there is no 4-arc-transitive graph of odd order.
Moreover, it was shown in \cite{Li1} that an $s$-arc-transitive graph of odd order with $s=2$ or 3 is
a normal cover of some $(G,2)$-arc-transitive graph where $G$ is an almost simple group,
led to the problem:

\ \ \ \ {\em Classify $(G,2)$-arc-transitive graphs of odd order with $G$ almost simple}.

This is one of a series of papers aiming to solve this problem, and does this work for alternating groups and symmetric groups.
The first one \cite{odd-excep} of the series of papers solves the problem for the exceptional groups of Lie type, and the sequel will solve the problem for other families of almost simple groups.

Let $\Ga=(V,E)$ be a connected $(G,2)$-arc-transitive graph of odd order, where $G$ is an almost simple group with socle being an alternating group. For the case where $G$ is primitive on $V$, it is easily deduced from \cite{P-W} that $\Ga$ is one of the complete graphs and the odd graphs. The main result of this paper shows that these   are all the graphs we  expected.
%The main result of this paper is the following theorem.

\begin{theorem}\label{zthm1}
Let  $G$ be an almost simple group with socle being an alternating group $\A_n$, and let $\Ga$ be a connected $(G,2)$-arc-transitive graph of odd order.
Then either
\begin{itemize}
\item[(i)] $\Ga$ is the complete graph $\K_n$, and $n$ is odd; or

\item[(ii)] $\Ga$ is the odd graph $\O_{2^e-1}$, and $n={2^{e+1}-1\choose 2^e-1}$ for some integer $e\geqslant 2$.

\end{itemize}
\end{theorem}
\begin{remark}
{\rm It would be infeasible to extend the classification in Theorem \ref{zthm1} to those graphs of  even order.
This is demonstrated by the work of Praeger-Wang in \cite{P-W} which presents a description of $(G,2)$-arc-transitive and $G$-vertex-primitive graphs with socle of $G$ being an alternating group.
}
\end{remark}

As a byproduct,
%Every finite simple group is contained in some alternating group,
%and alternating groups and symmetric groups undoubtedly bear a complicated subgroup structure.
%However,
the following result shows that subgroups of alternating groups and symmetric groups of odd index are very restricted: each insoluble composition factor is alternating except for three small exceptions.

\begin{theorem}\label{Alt-subgps}
Let  $G$ be an almost simple group with socle  $\A_n$,  and let $H$ be an insoluble proper subgroup of $G$ of odd index.
Then $G\in \{\A_n,\S_n\}$ and either
\begin{itemize}
\item[(i)] every insoluble composition factor of $H$ is an alternating group; or

\item[(ii)] $(G,H)=(\A_7,\GL(3,2))$, $(\A_8,\AGL(3,2))$ or $(\A_9,\AGL(3,2))$.

\end{itemize}
\end{theorem}

%See next section for the definition of {\it odd graphs}.
The notation used in the paper is standard, see for example the Atlas \cite{Atlas}.
In particular, a positive integer $n$ sometimes denotes a cyclic group of order $n$, and
for a prime $p$, the symbol $p^n$ sometimes denotes an elementary abelian $p$-group. For groups $A$ and $B$, an upward extension of $A$ by $B$ is denoted by $A{.}B$, and a semi-direct product of $A$ by $B$ is denoted by $A{:}B$.

For a positive integer $n$ and a prime $p$, let $n_p$ denote the $p$-part of $n$, that is, $n=n_pn'$ such that $n_p$ is a power of $p$ and $\gcd(n_p,n')=1$.
For a subgroup $H$ of a group $G$, let $|G:H|=|G|/|H|$, the index of $H$ in $G$, and denote by $\N_G(H)$ and $\C_G(H)$ the normalizer and the centralizer of $H$ in $G$, respectively.

\vskip 20pt

\section{Examples}\label{exam-sec}

We study the graphs which appear in our classification.

It is easily shown that, for an integer $n\geqslant  3$, the complete graph $\K_n$ is $(G,2)$-arc-transitive if and only if $G$ is a $3$-transitive permutation group of degree $n$.
Thus, if $n\geqslant  5$ is odd then $\K_n$ is one of the desired graphs.

The second type of example is the odd graph, defined below.

\begin{example}\label{exam-odd}
{\rm
Let $\Ome=\{1,2,\dots,2m+1\}$, and let $\Ome^{\{m\}}$ consist of $m$-subsets of $\Ome$.
Define a graph $(V,E)$ with vertex set and edge set
\[
V=\Ome^{\{m\}},\,
E=\{(\a,\b)\mid \a\cap \b=\emptyset\},
\]
respectively, which is called an {\it odd graph} and denoted by $\O_m$.}
\end{example}

The graph $\O_m$ has valency $m+1$, and has $\Sym(\Ome)=\S_{2m+1}$ to be the automorphism group, see \cite[pp. 147, Corollary 7.8.2]{Godsil-Royle}.
The order of $\O_m$   is given by
\[|V|=|\Ome^{\{m\}}|={2m+1\choose m}={(2m+1)!\over m!(m+1)!}.\]
For example, the Petersen graph is $\O_2$, which has order ${5\choose2}=10$ and valency 3;
$\O_3$ has order ${7\choose 3}=35$ and valency 4.
The former has even order, and the latter has odd order.
We next give a necessary and sufficient condition for ${2m+1\choose m}$ to be odd.

For a positive integer $n$, letting $2^{t+1}>n\geqslant 2^t$ for some integer $t\geqslant  0$, set
\[s(n)=\left[\frac{n}{2}\right]+\left[\frac{n}{2^2}\right]+\cdots+\left[\frac{n}{2^i}\right]+\cdots+\left[{n\over2^t}\right],\]
where $[x]$ is the largest integer which is not larger than $x$.
Then $[{n\over 2^i}]$ is the number of integers in $\{1,2,\dots,n\}$ which are divisible by $2^i$, and it follows that the 2-part of $n!$ is equal to $2^{s(n)}$.
Clearly, $2^{s(n)}=2^{s(n-1)}n_2$ if $n\geqslant  2$, where $n_2$ is the $2$-part of $n$.
We  observe that $[\frac{m}{2^i}]+[\frac{n}{2^i}]\leqslant [\frac{m+n}{2^i}]$ for all positive integers $i$.
It follows that
\begin{equation}\label{s(n)}
s(m)+s(n)\leqslant s(m+n),
\end{equation}
and
\begin{equation}\label{s(n)-1}
s(m)+s(n)=s(m+n) \Longleftrightarrow
\left[\frac{m}{2^i}\right]+\left[\frac{n}{2^i}\right]=\left[\frac{m+n}{2^i}\right] \mbox{ for all $i\geqslant 1$}.
\end{equation}
Further, if $s(m)+s(n)=s(m+n)$ then at least one of $n$ and $m$ is even.

Let $1\leqslant  m\leqslant  n$ and $\left[\frac{m}{2^i}\right]+\left[\frac{n}{2^i}\right]
=\left[\frac{m+n}{2^i}\right]$ for some $i\ge 1$.
Suppose that $a:=\left[\frac{m}{2^i}\right]\ne 0$. Then $b:=\left[\frac{n}{2^i}\right]\geqslant  a$.
Write $m=a2^i+c$ and  $n=b2^i+d$
for $c,d<2^i$.
We have  \[\left[\frac{m+n}{2^{i+1}}\right]=\left[{a+b\over 2}+{c+d\over 2^{i+1}}\right]\geqslant  \left[{a
+b\over 2}\right]\geqslant  \left[a\over 2\right]+\left[b\over 2\right]=\left[\frac{m}{2^{i+1}}\right]+\left[\frac{n}{2^{i+1}}\right].\]
Noting that $\left[{a+b\over 2}\right]\geqslant  1$, if   $\left[\frac{m+n}{2^{i+1}}\right]=\left[\frac{m}{2^{i+1}}\right]+\left[\frac{n}{2^{i+1}}\right]$
then $b\geqslant  2$, and so $\left[\frac{n}{2^{i+1}}\right]\ne 0$. Then, using (\ref{s(n)}) and (\ref{s(n)-1}), we have the following lemma.
\begin{lemma}\label{s(n)-2}
 Assume that
$s(m+n)=s(m)+s(n)$. If $m\leqslant  n$ and $\left[\frac{m}{2^i}\right]\ne 0$ then $\left[\frac{n}{2^{i+1}}\right]\ne 0$;
in particular, $m<n$, and $n\geqslant  2^t$ if  $\left[\frac{m+n}{2^{t}}\right]\ne 0$.
\end{lemma}

The following is a criterion for ${2m+1\choose m}$ to be odd.

\begin{lem}  \label{s(t)}
The number ${2m+1\choose m}={(2m+1)!\over m!(m+1)!}$ is odd if and only if $m+1$ is a $2$-power.
\end{lem}
\proof
Suppose that ${2m+1\choose m}$ is odd. Then $s(2m+1)=s(m)+s(m+1)$. Write $2^k\leqslant  m<2^{k+1}$. By Lemma
\ref{s(n)-2}, $\left[\frac{m+1}{2^{k+1}}\right]\ne 0$,
yielding $m+1\geqslant  2^{k+1}$, and so $m+1= 2^{k+1}$.

Conversely, we assume $m+1=2^\ell$ for some positive integer $\ell$.
Since $m=2^\ell-1$ and $2m+1=2^{\ell+1}-1$, we obtain
\[\left[{m\over2^i}\right]=\left[{2^\ell-1\over 2^i}\right]=\left\{
\begin{array}{ll}
2^{\ell-i}-1, &\mbox{for $1\leqslant i\leqslant \ell-1$,}\\
0, &\mbox{for $i\geqslant \ell$.}
\end{array}\right.\]
\[\left[{2m+1\over2^i}\right]=\left[{2^{\ell+1}-1\over 2^i}\right]=\left\{
\begin{array}{ll}
2^{\ell+1-i}-1, &\mbox{for $1\leqslant i\leqslant \ell$,}\\
0, &\mbox{for $i\geqslant \ell+1$.}
\end{array}\right.\]
Therefore, we have
\[\begin{array}{rll}
s(m)&=&(2^{\ell-1}-1)+(2^{\ell-2}-1)+\dots+(2-1),\\

s(m+1)&=&2^{\ell-1}+2^{\ell-2}+\dots+2+1,\\

s(2m+1)&=&(2^{\ell+1-1}-1)+(2^{\ell+1-2}-1)+\dots+(2-1).\\

\end{array}\]
Then $s(m)+s(m+1)=s(2m+1)$, and ${2m+1\choose m}$ is odd.
\qed

By the above lemma, we   get the following consequence.

\begin{corollary}\label{symmetric index odd}
The odd graph $\O_m$ is of odd order if and only if $m+1$ is a $2$-power.
\end{corollary}

\vskip 30pt

\section{Subgroups with odd index in $\A_n$ or $\S_n$}\label{proof-th1}

Let  $G$ be an almost simple group with socle  $\A_n$. Then either $G\in \{\A_n,\S_n\}$ or $n=6$ and $G\in\{\PGL(2,9),\M_{10},\PGammaL(2,9)\}$.
In this section, we shall determine the insoluble composition factors of subgroups of $G$ of odd index.

For the natural action of $\S_n$ on $\Ome=\{1,2, \ldots, n\}$ and a subset $\Del\subseteq \Ome$,
the symmetric group $\Sym(\Del)$ is sometimes identified with a subgroup of $\S_n$.
Thus we write the set-stabilizer $G_\Del$ as $(\Sym(\Del)\times \Sym(\Ome\setminus \Del))\cap G$ or simply, $G_\Del=(\S_m\times \S_{n-m})\cap G$ if $|\Del|=m$.
Also, $(\S_{m}\wr \S_k )\cap G$ stands for the  stabilizer in $G$ of some
 partition  of $\Ome$ into $k$ parts  with equal size $m$.

Based on O'Nan-Scott theorem, the following lemma was first obtained by Liebeck and Saxl \cite{LS}.

\begin{lemma}[\cite{LS})]\label{max-subgps}
Let $G$ have socle $T=\A_n$ with $n\geqslant 5$ and have a maximal subgroup $M$ of odd index.
Then one of the  following holds:
\begin{itemize}
\item[(1)] $M=(\S_m\times\S_{n-m})\cap G$ with $1\leqslant  m<{n\over 2}$; or

\item[(2)] $M=(\S_m\wr\S_k)\cap G$, where  $n=mk$ and $m,k>1$; or

\item[(3)] $G=\A_7$ and $M\cong \SL(3,2)$, or $G=\A_8$ and $M\cong \AGL(3,2)$; or

\item[(4)] $G=\PGL(2,9)$, $\M_{10}$ or $\PGammaL(2,9)$, and $M$ is a Sylow
$2$-subgroup of $G$.
\end{itemize}
In particular, if $G\ne \A_7$ or $\A_8$, then each insoluble composition factor of $M$ is an alternating group.
\end{lemma}

For a subgroup $X\leqslant  \S_n$ fixing a subset $\Del\subseteq\Ome$, denote by $X^\Del$ the permutation group induced by $X$ on $\Del$.

\begin{lemma}\label{tech}
Let $G=\S_n$ or $\A_n$ with $n\geqslant 5$, and let $H$ be a subgroup of $G$ with odd index $|G:H|>1$.
Suppose that  $H$ normalizes  a subgroup $L=\Sym({\Del_1})\times \cdots\times\Sym({\Del_t})$ of $\S_n$, where $t\geqslant  2$ and
$\Ome=\cup_{i=1}^t\Del_i$. Then
\begin{enumerate}
\item[(1)] $|(L\cap G):(L\cap H)|$ and $|(L\cap G)^{\Del_i}:(L\cap H)^{\Del_i}|$ are odd, where
$1\leqslant  i\leqslant  t$;

\item[(2)] each composition factor   of $L\cap H$ is a composition factor of
 some $(L\cap H)^{\Del_i}$.
\end{enumerate}
\end{lemma}
\proof
By the assumption $LH$ is a subgroup of $\S_n$, and so $H\leqslant  LH\cap G=(L\cap G)H\leqslant  G$.
Thus $|(L\cap G)H:H|$ is odd. Then $|(L\cap G):(L\cap H)|$ is odd as $|(L\cap G)H:H|={|L\cap G|\over |L\cap H|}$.

Let $L_i$ be the kernel of $L\cap G$ acting on $\Del_i$, where
$1\leqslant  i\leqslant  t$. Then $L^{\Del_i}\cong L/L_i$, $(L\cap G)^{\Del_i}\cong (L\cap G)/(L_i\cap G)$ and $(L\cap H)^{\Del_i}\cong (L\cap H)(L_i\cap G)/(L_i\cap G)$.
Since $|(L\cap G):(L\cap H)|$ is odd, $|(L\cap G):(L\cap H)(L_i\cap G)|$ is odd, and so is $|(L\cap G)^{\Del_i}:(L\cap H)^{\Del_i}|$, as in part~(1).

Let $S$ be a composition factor of $L\cap H$. Since $(L\cap H)^{\Del_t}\cong (L\cap H)(L_t\cap G)/(L_t\cap G)\cong (L\cap H)/ (L_t\cap H)$, it follows that
 $S$  is a composition factor of one of $(L\cap H)^{\Del_t}$ and $L_t\cap H$. If $S$  is a composition factor of  $(L\cap H)^{\Del_t}$, then part (2) holds by taking $i=t$.
Now let $S$  be a composition factor of   $L_t\cap H$, and consider the triple $(L_t, L_t\cap G, L_t\cap H)$.
By induction, we may assume that $S$ is a composition factor of $(L_t\cap H)^{\Del_i}$ for some $i\leqslant  t-1$.
Since $L_t\cap H\unlhd L\cap H$, we have $(L_t\cap H)^{\Del_i}\unlhd (L\cap H)^{\Del_i}$, and
thus $S$ is a composition factor of $(L\cap H)^{\Del_i}$. Then part (2) follows.
\qed

Now we prove Theorem \ref{Alt-subgps} for $G=\S_n$.

\begin{lemma}\label{Sym}
Let $G=\S_n$ with $n\geqslant 5$, and let $H$ be an insoluble subgroup of $G$ with odd index $|G:H|>1$.
Then each insoluble composition factor of $H$ is an alternating group.
\end{lemma}
\proof
We prove this lemma by induction on $n$.
Let $S$ be an insoluble composition factor of $H$.
Take a maximal subgroup $M$ of $G$ with $H\leqslant M$.
By Lemma~\ref{max-subgps}, either $M=\S_m\times\S_{n-m}$ with $1\leqslant m<n/2$, or $M=\S_m\wr\S_k$ with $mk=n$ and $m,k>1$.

For $M=\S_m\times\S_{n-m}$,   Lemma \ref{tech} works for $H$ and $M$,
which yields that $S$ is a composition factor of a subgroup with odd index in $\S_k$ for some $k<n$,
and the lemma holds by induction.
Thus, let $M=\S_m\wr\S_k$ with $mk=n$ and $m,k>1$ in the following.

Let $L$ be the base subgroup of the wreath product $\S_m\wr\S_k$. Then Lemma \ref{tech}
works for the triple $(L,H, L\cap H)$, and hence the lemma holds by induction if $S$ is  a composition factor of
$L\cap H$.

Assume that $S$ is not a composition factor of
$L\cap H$. Then $S$ is a composition factor of $H/(L\cap H)$. Noting that $HL/L\cong H/(L\cap H)$, it implies that
$S$ is  a composition factor of   $HL/L$. Consider that pair $M/L$ and $HL/L$.
Since $|G:H|$ is odd, $|M:(HL)|$ and hence $|(M/L):(HL/L)|$ is also odd. Further, $M/L\cong \S_k$. Then, since $k<n$, the lemma holds by induction.
\qed

Now we handle the case $G=\A_n$.

\begin{lemma}\label{Alt}
Let $G=\A_n$ with $n\geqslant  5$.
Let $H$ be an insoluble  subgroup of $G$ with odd index $|G:H|>1$.
Then either
\begin{itemize}
\item[(i)] $(G,H)$ is one of $(\A_7,\GL(3,2))$, $(\A_8,\AGL(3,2))$ and $(\A_9,\AGL(3,2))$; or
\item[(ii)] every insoluble composition factor of $H$ is an alternating group.
\end{itemize}
\end{lemma}
\proof
If $n\leqslant  9$ then the lemma is easily shown by checking the subgroups of $\A_n$.
In the following, by induction on $n$, we show (ii) of this lemma always holds for $n\geqslant  10$.

Let $n\geqslant  10$, and let $S$ be an   insoluble composition factor of $H$.
Take a maximal subgroup $M$ of $\A_n$ with $H\leqslant  M$. By Lemma \ref{max-subgps}, $M=(\S_m\times\S_{n-m})\cap \A_n$ with $1\leqslant  m<n/2$, or $M=(\S_m\wr\S_k)\cap \A_n$ with $mk=n$ and $m,k>1$.

Suppose that $n=10$. Then $M\cong \S_8$ or $2^4{:}\S_5$.
By the Atlas \cite{Atlas}, $\S_8$ has no insoluble subgroup of odd index. Then  $M\cong 2^4{:}\S_5$, and we have $S=\A_5$.
Thus, in the following, we  let  $n\geqslant  11$, and process in two cases.

{\bf Case 1}. Let $M=(\S_m\times\S_{n-m})\cap \A_n$.
If $m=1$ then $M=\A_{n-1}$ and, since $10\leqslant  n-1<n$, $S$ is alternating by induction.
Now let
 $m\geqslant  2$. Writing $M=(\Sym(\Del)\times \Sym(\Ome\setminus\Del))\cap \A_n$ with $|\Del|=m$,
we have $M=(\Alt(\Del)\times \Alt(\Ome\setminus\Del))\l \s_1\s_2\r$, where $\s_1\in \Sym(\Del)$ and
$\s_2\in \Sym(\Ome\setminus \Del)$ are transpositions. Then $M^\Del\cong \S_m$ and $M^{\Ome\setminus \Del}\cong \S_{n-m}$.
By  Lemma \ref{tech},
  $S$ is a composition factor of a subgroup with odd index in either $\S_m$  or  $\S_{n-m}$.  Then $S$ is alternating by Lemma \ref{Sym}.

{\bf Case 2}.
Let  $M=(\S_m\wr\S_k)\cap \A_n$. Let $L=\S_m^k$ be the base group of the wreath product $\S_m\wr\S_k$.
Note that $S$ is a  composition factor of one of $H/(L\cap H)$ and $L\cap H$.

Assume that $S$  is a  composition factor of $H/(L\cap H)$.
Then $S$ is  a composition factor of   $HL/L$ as
 $HL/L\cong H/(L\cap H)$.
It is easily shown that $|(M/L):(HL/L)|$ is  odd. Further, since $M/L\cong \S_k$, we know that $S$ is alternating by Lemma \ref{Sym}.

Now let $S$  be a  composition factor of $L\cap H$.
Write $L=\Sym(\Del_1)\times\cdots\times \Sym(\Del_k)$, where $|\Del_i|=m$.
Then $L\cap \A_n=(\Alt(\Del_1)\times\cdots\times \Alt(\Del_k))\l\s_1\s_t,\s_2\s_t,\ldots,\s_{t-1}\s_t\r$, where $\s_i\in \Alt(\Del_i)$ are transpositions. It follows that $(L\cap \A_n)^{\Del_i}\cong \S_m$ for $1\leqslant  i\leqslant  t$.
Thus, using Lemmas \ref{tech} and \ref{Sym}, $S$ is an alternating group.
\qed

Finally, if $n=6$ and $G=\PGL(2,9)$, $\M_{10}$ or $\PGammaL(2,9)$ then, by Lemma \ref{max-subgps}, $G$ has no insoluble proper subgroup of odd index.
The proof of Theorem~\ref{Alt-subgps} now follows from Lemmas~\ref{Sym} and \ref{Alt}.

\vskip 20pt

\section{$2$-Arc-transitive graphs}\label{sect=proof-th3}

In this section, we assume that $\Ga=(V,E)$ is a connected $(G,2)$-arc-transitive graph of odd order and valency at least $3$, where $G\leqslant \Aut\Ga$.

\subsection{Stabilizers}\label{sub=stab}
Fix a 2-arc $(\a,\b,\g)$ of $\Ga$.
Let $G_\a$ be the stabilizer of $\a$ in $G$.
Then $G_\a$ acts 2-transitively on the neighborhood $\Ga(\a)$ of $\a$ in $\Ga$.
Let $G_\a^{[1]}$ be the kernel of $G_\a$ on $\Ga(\a)$, and let $G_\a^{\Ga(\a)}$ be the 2-transitive permutation group induced by $G_\a$ on $\Ga(\a)$.
Then $G_\a^{\Ga(\a)}\cong G_\a/G_\a^{[1]}$.
Clearly, $G_\a^{[1]}\unlhd  G_{\a\b}$, and
\begin{equation}\label{exten}
(G_\a^{[1]})^{\Ga(\b)}\unlhd  G_{\a\b}^{\Ga(\b)}\cong G_{\a\b}^{\Ga(\a)}.
\end{equation}
Let $G_{\a\b}^{[1]}=G_\a^{[1]}\cap G_\b^{[1]}$, the point-wise stabilizer of the `double star' $\Ga(\a)\cup\Ga(\b)$.
A fundamental result about 2-arc-transitive graphs characterizes $G_{\a\b}^{[1]}$.

\begin{theorem}\label{double-star}{\rm (Thompson-Wielandt Theorem)}
$G_{\a\b}^{[1]}$ is a $p$-group with $p$ prime.
\end{theorem}

By definition, we have $G_{\a\b}^{[1]}\unlhd  G_\b^{[1]}\unlhd  G_{\b\g}$, and so
\[(G_{\a\b}^{[1]})^{\Ga(\g)}\unlhd  (G_\b^{[1]})^{\Ga(\g)}\unlhd  G_{\b\g}^{\Ga(\g)}.\]
Let $O_p((G_\b^{[1]})^{\Ga(\g)})$ and $O_p(G_{\b\g}^{\Ga(\g)})$ be the maximal normal $p$-subgroups of $(G_\b^{[1]})^{\Ga(\g)}$ and $G_{\b\g}^{\Ga(\g)}$, respectively. Then
\[(G_{\a\b}^{[1]})^{\Ga(\g)}\unlhd  O_p((G_\b^{[1]})^{\Ga(\g)})\unlhd  O_p(G_{\b\g}^{\Ga(\g)}).\]
Suppose that $(G_{\a\b}^{[1]})^{\Ga(\g)}=1$. Then $G_{\a\b}^{[1]}\leqslant  G_\g^{[1]}$, and so $G_{\a\b}^{[1]}\leqslant  G_{\b\g}^{[1]}$. Noting that
$G_{\a\b}^{[1]}\cong  G_{\b\g}^{[1]}$, we have $G_{\a\b}^{[1]}= G_{\b\g}^{[1]}$. Then the connectedness of $\Ga$ yields that $G_{\a\b}^{[1]}=G_{\a'\b'}^{[1]}$ for each arc $(\a',\b')$ of $\Ga$, and hence $G_{\a\b}^{[1]}=1$.
Thus, if $G_{\a\b}^{[1]}$ is a non-trivial $p$-group, then so is $(G_{\a\b}^{[1]})^{\Ga(\g)}$, and then $O_p(G_{\b\g}^{\Ga(\g)})\ne 1$.
Noting that $G_{\a\b}^{\Ga(\a)}\cong G_{\b\g}^{\Ga(\g)}$, we have a useful conclusion.

\begin{lemma}\label{p-local}
Let $\{\a,\b\}\in E$.
If $G_{\a\b}^{[1]}$ is a nontrivial  $p$-subgroup, then $G_{\a\b}^{\Ga(\a)}$ has a nontrivial normal $p$-subgroup, where $p$ is a prime.
\end{lemma}

Recall that $G_\a^{\Ga(\a)}$ is  2-transitive on $\Ga(\a)$.
Inspecting 2-transitive permutation groups (refer to \cite[page 194-197, Tables 7.3 and 7.4]{Cameron}), we have the following result.

\begin{lemma}\label{stab-1}
Let $G$ be an almost simple group with socle   $\A_n$,
and $\{\a,\b\}\in E$.
Then either $G_\a$ is soluble, or $G\in\{\A_n,\S_n\}$  and one of the following holds.
\begin{itemize}
\item[(1)] $\soc(G_\a^{\Ga(\a)})\cong \A_m$ for some $m\geqslant 5$, and one of the following holds:
\begin{enumerate}
\item[(i)] $G_\a^{\Ga(\a)}\cong \A_m$ or $\S_m$ for even $m\geqslant  6$, and $G_{\a\b}^{\Ga(\a)}\cong \A_{m-1}$ or $\S_{m-1}$, respectively;
\item[(ii)] $G_\a^{\Ga(\a)}\cong \PSL(2,5)$ or $\PGL(2,5)$, and $G_{\a\b}^{\Ga(\a)}\cong \D_{10}$ or $5{:}4$, respectively;
\item[(iii)] $G_\a^{\Ga(\a)}\cong \PSL(2,9).\calO$, and $G_{\a\b}^{\Ga(\a)}\cong 3^2{:}(4.\calO)$, where $\calO\leqslant 2^2$.
\end{enumerate}

\item[(2)] $G_\a^{\Ga(\a)}\cong 2^4{:}H$, where $H=G_{\a\b}^{\Ga(\a)}\cong \A_5$, $\S_5$, $3\times \A_5$, $(3\times \A_5).2$, $\A_6$, $\S_6$, $\A_7$ or $\A_8$; in particular, $G_{\a\b}^{[1]}=1$.

\end{itemize}
\end{lemma}
\proof
Note that \begin{equation}\label{extension} G_\a=G_\a^{[1]}.G_\a^{\Ga(\a)}=(G_{\a\b}^{[1]}.(G_\a^{[1]})^{\Ga(\b)}).G_\a^{\Ga(\a)}.
\end{equation}
Clearly, if $G_\a^{\Ga(\a)}$ is insoluble then $G_\a$ is insoluble.
If $G_\a^{\Ga(\a)}$ is soluble then, by (\ref{exten}), $(G_\a^{[1]})^{\Ga(\b)}$ is soluble, and so $G_\a$ is soluble by (\ref{extension}).
Thus   $G_\a$ is soluble if and only if   $G_\a^{\Ga(\a)}$ is soluble.
To finish the proof of this lemma, we assume that  $G_\a$ is insoluble in the following; in particular,  $G\in\{\A_n,\S_n\}$ by Theorem \ref{Alt-subgps}.
Since $\Ga$ is $(G,2)$-arc-transitive,
$G_\a^{\Ga(\a)}$ is an insoluble $2$-transitive permutation group.
As $|V|$ is odd, the valency $|\Ga(\a)|$ is even, and so $G_\a^{\Ga(\a)}$ is of even degree.

{\bf Case 1}. First assume that $G_\a^{\Ga(\a)}$ is an almost simple $2$-transitive permutation group with socle $S$ say.
By Theorem~\ref{Alt-subgps},  either  $S\cong \A_m$ for some $m\geqslant  5$, or one of the following cases occurs:
\begin{enumerate}
\item[(a)] $G=\A_7$, $G_\a=\SL(3,2)$;
\item[(b)] $G=\A_8$,  $G_\a=\AGL(3,2)$;
\item[(c)] $G=\A_9$,  $G_\a=\AGL(3,2)$.
\end{enumerate}

For (a) and (b), we have that $|V|=15$, and $G$ is $2$-transitive on $V$, yielding $\Ga\cong \K_{15}$.
Noting that $\Ga$ is $(G,2)$-arc-transitive, it follows that $G=\A_7$ or $\A_8$ is $3$-transitive on the $15$ vertices of $\Ga$, which is impossible.

Suppose that (c) occurs.
Let $G_\a^{\Ga(\a)}$ be of affine type. Then $G_{\a\b}=\SL(3,2)$; in this case, as a subgroup, $\SL(3,2)$ is self-normalized in $\A_9$.
Thus there is no element in $G$ interchanging $\a$ and $\b$, which contradicts the arc-transitivity of $G$ on $\Ga$.
Thus $G_\a^{\Ga(\a)}$ is almost simple. Then $G_\a^{[1]}=\ZZ_2^3$ and $G_\a^{\Ga(\a)}\cong \SL(3,2)\cong \PSL(2,7)$. Since $\Ga$ has even valency, considering the $2$-transitive permutation representations of $\SL(3,2)$, we have $|\Ga(\a)|=8$.
Then $G_\a^{[1]}$ is not faithful on $\Ga(\b)\setminus\{\a\}$, and so $G_{\a\b}^{[1]}$ is a non-trivial normal $2$-group. By Lemma \ref{double-star},
$G_{\a\b}^{\Ga(\a)}$ has a non-trivial $2$-subgroup; however,
$G_{\a\b}^{\Ga(\a)}\cong \ZZ_7{:}\ZZ_3$, a contradiction.

Let $S\cong \A_m$.
Note that $\A_5\cong \PSL(2,5)$ and $\A_6\cong \PSL(2,9)$.
By the classification of 2-transitive permutation groups (refer to \cite[page 197, Table 7.4]{Cameron}), since $|\Ga(\a)|$ is even, either $|\Ga(\a)|=m$ with $m$ even,
or $(S,|\Ga(\a)|)$ is one of $(\PSL(2,5),6)$ and $(\PSL(2,9),10)$.
Then part~(1) follows.

{\bf Case 2}. Now suppose that $G_\a^{\Ga(\a)}$ is an insoluble affine group.
Then $|\Ga(\a)|=2^d$ for some positive integer $d\geqslant  3$, and $G_{\a\b}^{\Ga(\a)}\leqslant  \GL(d,2)$. In particular, by \cite{Weiss}, we have $G_{\a\b}^{[1]}=1$.
Since each insoluble composition factor of $G_\a^{\Ga(\a)}$ is alternating,
by the classification of affine 2-transitive permutation groups (see \cite[page 195, Table 7.3]{Cameron}), we conclude that $d=4$ and
$G_{\a\b}^{\Ga(\a)}$ is isomorphic to one of $\A_5$ (isomorphic to
$\SL(2,4)$), $\S_5$ (isomorphic to
$\SigmaL(2,4)$), $\ZZ_3\times \A_5$ (isomorphic to
$\GL(2,4)$), $(\ZZ_3\times \A_5).2$ (isomorphic to
$\GammaL(2,4)$),  $\A_6$ (isomorphic to $\Sp(4,2)'$),  $\S_6$ (isomorphic to $\Sp(4,2)$), $\A_7$ and $\A_8$
(isomorphic to $\GL(4,2)$).
This gives rise to the candidates in part~(2).
\qed

Let $G$ be an almost simple group with socle  $\A_n$. We next organize our analysis of the candidates for $G_\a$ according to the description in Lemma~\ref{stab-1}.
Note that $G\in\{\A_n,\S_n\}$ if $G_\a$ is insoluble.

\subsection{Almost simple stabilizers}
Assume that $G_\a^{\Ga(\a)}$ is almost simple, where $\a\in V$.
First we consider the candidates in Lemma~\ref{stab-1}\,(1)(i).

\begin{lemma}\label{main-case}
Let $\{\a,\b\}\in E$.
Assume $G_\a^{\Ga(\a)}\cong \A_m$ or $\S_m$, and $G_{\a\b}^{\Ga(\a)}\cong \A_{m-1}$ or $\S_{m-1}$, respectively, where $|\Ga(\a)|=m\geqslant6$ is even.
Then one of the following holds:
\begin{itemize}

\item[(i)] $(G_\a,G)=(\A_m,\A_{m+1})$ or $(\S_m,\S_{m+1})$, and $\Ga=\K_{m+1}$, where $m$ is even;

\item[(ii)] $G_\a=(\S_m\times\S_{m-1})\cap G$, $G=\A_{2m-1}$ or $\S_{2m-1}$, respectively, and $\Ga=\O_{m-1}$, where $m$ is a power of $2$.

\end{itemize}
\end{lemma}
\proof
Since $G_{\a\b}^{\Ga(\a)}$ is almost simple, $G_{\a\b}^{[1]}=1$ by Lemma~\ref{p-local}, and so
\begin{equation}\label{double-star=1} G_\a=G_\a^{[1]}.G_\a^{\Ga(\a)}=(G_{\a\b}^{[1]}.(G_\a^{[1]})^{\Ga(\b)}).G_\a^{\Ga(\a)}=(G_\a^{[1]})^{\Ga(\b)}.G_\a^{\Ga(\a)}.
\end{equation}
Since $(G_\a^{[1]})^{\Ga(\b)}$ is isomorphic to a normal subgroup of $G_{\a\b}^{\Ga(\a)}$, we have $(G_\a^{[1]})^{\Ga(\b)}=1$, or $(G_\a^{[1]})^{\Ga(\b)}\cong \A_{m-1}$ or $\S_{m-1}$.
It follows that $G_\a\cong \A_m$, $\S_m$, $\A_{m-1}\times \A_m$,  $(\A_{m-1}\times \A_m).2$ or  $\S_{m-1}\times \S_m$.

{\bf Case 1.} Assume first that $G_\a\cong \A_m$ or $\S_m$, where $m$ is even.
Since $G=\A_n$ or $\S_n$ and $|G:G_\a|$ is odd, it follows that
either $n=m+1$ and $G_\a=\S_m\cap G$, or $n=m+k$, $G=\A_{m+k}$ and $G_\a\cong \S_m$ for $k\in \{2,3\}$.

Suppose that $n=m+k$, $G=\A_{m+k}$  and $G_\a\cong \S_m$, where $k=2$ or $3$.
Then $G_{\a\b}\cong \S_{m-1}$ since $\Ga$ is of valency $m$.
Consider the maximal subgroups of $G=\A_{m+k}$ which contains $G_\a$.
By Lemma \ref{max-subgps}, we conclude that $G_\a$ is contained in the stabilizer of an $m$-subset of $\Ome=\{1,2,\dots,m+k\}$, say $\Del=\{1,2,\dots,m\}$.
Thus we may let
$G_\a=\Alt(\Del).\l\s\r$, where $\s=(1\,\,2)(m+1\,\,m+k)$.
Without loss of generality, we may assume that $G_{\a\b}=\Alt(\Del\setminus\{m\}).\l\s\r$.
Let $g\in G$ interchange $\a$ and $\b$.
Then $g$ normalizes $G_{\a\b}$, and hence $g$ fixes $\Del\setminus\{m\}$ setwise, and $\s^g=(i\,\,j)(m+1\,\,m+k)$.
It follows that $\Del$ and $\{m+1,m+k\}$ are two orbits of $\l G_\a,g\r$, which is a contradiction since $\l G_\a,g\r$ should be equal to $G$.
Thus $(G_\a,G)=(\A_m,\A_{m+1})$ or $(\S_m,\S_{m+1})$.
It then follows that $\Ga=\K_{m+1}$, as in part~(i).

{\bf Case 2.} Now assume that $G_\a$ has a subgroup isomorphic to $\A_m\times\A_{m-1}$.
Clearly, $n\geqslant  2m-1$.
Recall that $2^{s(l)}$ is the 2-part of $l!$, see Section \ref{exam-sec}.
Then $|G|_2\geqslant  2^{s(n)-1}$ and $|G_\a|_2\leqslant  2^{s(m)+s(m-1)}$.
Since $|G:G_\a|$ is odd, $s(m)+s(m-1)\geqslant  s(n)-1\geqslant  s(2m-1)-1$.
By (\ref{s(n)}) given in Section \ref{exam-sec},
$s(2m-1)\geqslant  s(m)+s(m-1)$, and so
\[s(2m-1)\geqslant  s(m)+s(m-1)\geqslant  s(n)-1\geqslant  s(2m-1)-1.\]
Since $m$ is even, $2m$ is divisible by $2^2$, and hence $s(2m)\geqslant  s(2m-1)+2$.
It follows that $n<2m$. Therefore, we have
\[n=2m-1\]
 and $s(2m-1)=s(m)+s(m-1)$. Then $m$ is a power of $2$ by Lemma \ref{s(t)}.
Since $|G:G_\a|$ is odd, either $G=\A_{2m-1}$ and $G_\a=(\A_m\times\A_{m-1}).2$, or
$G=\S_{2m-1}$ and $G_\a=\S_m\times\S_{m-1}$.
That is to say, $G_\a$ is the stabilizer of $G$ acting on the set of $(m-1)$-subsets of $\{1,2,\dots,2m-1\}$.
It follows since $\Ga$ is $(G,2)$-arc-transitive that $\Ga=\O_{m-1}$ is an odd graph, as in part~(ii).
\qed

Next, we handle the candidates in part~(1)(ii-iii)   of Lemma~\ref{stab-1}.

\begin{lemma}\label{PSL(2,5)}
There is no $2$-arc-transitive graph corresponding to part~{\rm (1)(ii)}   of Lemma~{\rm \ref{stab-1}}.
\end{lemma}
\proof
Suppose that $G_\a^{\Ga(\a)}\cong \PSL(2,5)$ or $\PGL(2,5)$, and $G_{\a\b}^{\Ga(\a)}\cong \D_{10}$ or $5{:}4$.
By Lemma~\ref{p-local}, $G_{\a\b}^{[1]}$ is a 5-group, and so $|G_\a^{[1]}|_2=|(G_\a^{[1]})^{\Ga(\b)}|_2$ divides $|G_{\a\b}^{\Ga(\b)}|_2$.
Thus
\[|G_\a|_2=|G_\a^{[1]}|_2|G_\a^{\Ga(\a)}|_2\leqslant 2^5,\]
that is, a Sylow $2$-subgroup of $G_\a$ has  order a divisor of $2^5$.
It follows that $G\leqslant\S_7$.
Since $G_\a^{\Ga(\a)}\cong \PSL(2,5)$ or $\PGL(2,5)$, we conclude that either $G=\A_7$ and $G_\a\cong \S_5$, or $G=\S_7$ and $G_\a=\S_2\times\S_5$.
Then $\Ga$ is an orbital graph of $G=\S_7$ acting on 2-subsets of $\{1,2,\dots,7\}$, which is not 2-arc-transitive.
\qed

%The next lemma deals with the candidates in part~(1)(iii) of Lemma~\ref{stab-1}.

\begin{lemma}\label{PSL(2,9)}
There is no $2$-arc-transitive graph corresponding to to part~{\rm (1)(iii)}   of Lemma~{\rm \ref{stab-1}}.
\end{lemma}
\proof
Suppose that $G_\a^{\Ga(\a)}\cong \PSL(2,9).\calO$, and $G_{\a\b}^{\Ga(\a)}\cong 3^2{:}(4.\calO)$, where $\calO\leqslant 2^2$.
By Lemma~\ref{p-local}, $G_{\a\b}^{[1]}$ is a 3-group, and so $|G_\a^{[1]}|_2=|(G_\a^{[1]})^{\Ga(\b)}|_2$ divides $|G_{\a\b}^{\Ga(\b)}|_2$.
We have
\[|G_\a|_2=|G_\a^{[1]}|_2|G_\a^{\Ga(\a)}|_2\leqslant 2^9,\]
that is, a Sylow $2$-subgroup of $G_\a$ is of order dividing $2^9$.
It follows that $G\leqslant\A_{13}$, and further, either $G\leqslant\S_{11}$, or $G$ is one of $\A_{12}$ and $\A_{13}$.

Suppose $|G|_2=2^9$.
Then $G=\S_{11}$, $\A_{12}$ or $\A_{13}$, and moreover,
$G_\a^{\Ga(\a)}\cong \PSL(2,9).2^2$ and $G_\a^{[1]}\cong 3^2{:}[2^4]$, and hence
\[G_\a=(\PSL(2,9)\times (3^2{:}4)).[2^4].\]
By the Atlas \cite{Atlas}, $G$ does not have a subgroup of odd index which contains a normal subgroup $\PSL(2,9)\times (3^2{:}4)$, which is a contradiction.
Thus $|G|_2\leqslant 2^8$, and then $G\leqslant\A_{11}$ or $\S_{10}$.
Checking the subgroups of $G$ with odd index, we conclude that $\A_7\leqslant  G\leqslant  \S_7$ and $\A_6\leqslant  G_\a\leqslant  \S_6$.
It follows that $\Ga=\K_7$, which is not possible since $\Ga$ should have valency 10.
\qed

\subsection{The affine stabilizers} Let $\{\a,\b\}\in E$. Assume that $G_\a^{\Ga(\a)}$ is an affine $2$-transitive permutation group.

Now consider the case where $G_\a$ is soluble.
By \cite{odd-excep}, Theorem \ref{zthm1} holds for the case where $G_\a$ is soluble.

\begin{lemma}\label{soluble-case}
If $G_\a$ is soluble, then $\Ga$ has valency $4$, and either
\begin{itemize}
\item[(i)] $n=5$ and $\Ga$ is the complete graph $\K_5$, or
\item[(ii)]
$n=7$ and $\Ga$ is the odd graph $\O_3$ of order $35$.
\end{itemize}
\end{lemma}

We now consider   the candidates for $G_\a^{\Ga(\a)}$ in part~(2) of Lemma~\ref{stab-1}.

\begin{lemma}\label{stab-2}
There is no $2$-arc-transitive graph corresponding to  part~{\rm {(2)}}   of Lemma~{\rm \ref{stab-1}}.
\end{lemma}
\proof
Suppose that $G_\a^{\Ga(\a)}\cong 2^4{:}H$ is affine and described as in  part~(2) of Lemma~\ref{stab-1}.
Let $\{\a,\b\}\in E$.
Since $G_{\a\b}^{[1]}=1$, (\ref{exten}) yields that $G_\a^{[1]}$ is isomorphic to a normal subgroup of
$H=G_{\a\b}^{\Ga(\a)}$. Then the outer automorphism group of $G_\a^{[1]}$ has order at most $4$.
It follows that $G_\a$ has a (minimal) normal subgroup $N$  which is regular on $\Ga(\a)$, and thus
\[G_\a=N{:}G_{\a\b},\, \C_{G_\a}(N)=N\times G_\a^{[1]}.\]
Moreover, $|G_\a^{[1]}|_2$ is a  divisor of $|G_{\a\b}^{\Ga(\b)}|_2=|H|_2$,
and then $|G|_2=|G_\a|_2$ is a divisor of $2^4|H|_2^2$. In particular, $2^6\leqslant  |G|_2\leqslant  2^{16}$, and then $8\leqslant  n\leqslant 19$.

Consider the natural action  of $G_\a$ on $\Ome=\{1,2,\dots,n\}$, and choose a $G_\a$-orbit $\Del$ such that
$N$ is nontrivial on $\Del$. Let $|\Del|=m$. Then $m$ is even, and  $|G_\a^{\Del}|_2=|\S_m|_2$ or $|\A_m|_2$ by Lemma \ref{tech}.

Let $K$ be the kernel of
$G_\a$ acting on $\Del$. Then $K\cap N=1$ as $N$ is a minimal normal subgroup of $G_\a$, and so $K\le \C_{G_\a}(N)=N\times G_\a^{[1]}$.
It follows that $K\le G_\a^{[1]}$, and hence $G_\a^{\Del}$ is insoluble. In particular,  $m\geqslant 6$.

{\bf Case 1.} Suppose that $K$ is soluble. Then $|K|_2=1$, and  $2^4|H|_2||G_\a^{[1]}|_2=|G_\a|_2=|G_\a^{\Del}|_2=|\S_m|_2$ or $|\A_m|_2$.
Recalling that $|G_\a|_2=|G|_2=|\S_n|_2$ or $|\A_n|_2$, we have $n\leqslant m+3$.
If $N$ is transitive on $\Del$, then $m=|N|=16$, yielding $|G_\a|_2=2^{15}$ or $2^{14}$, which is impossible.
Thus $N$ is intransitive on $\Del$, and then $G_\a^{\Del}\lesssim\S_\ell\wr\S_k$, where $\ell,k>1$, $m=\ell k$ and $\ell$ is the size of each $N$-orbit. In particular, $\ell=2$, $4$ or $8$.

For $\ell=4$ or $8$, since $m=\ell k\leqslant n\leqslant 19$, we have $m=16$, which yields a contradiction as
above. Therefore, $\ell=2$ and, since $G_\a^{\Del}$ is insoluble, $5\leqslant k\leqslant 9$. Then
$G_\a$ has exactly one insoluble composition factor, and thus $|G_\a|_2=|G_\a^{\Del}|_2=2^4|H|_2$. This implies that $k=5$, $m=10$, and
$|G_\a|_2=2^7$ or $2^8$. Then $G=\A_{11}$ or $\A_{10}$, and $G_\a=2^4{:}\S_5$ which is faithful on $\Del$.
Thus $G_{\a\b}\cong \S_5$, which has two orbits on $\Del$ of equal size $5$.

Let $g\in G$ with $(\a,\b)^g=(\b,\a)$. Then $g$ normalizes $G_{\a\b}$, fixes $\Ome\setminus\Del$ and either interchanges or fixes those two  $G_{\a\b}$-orbits on $\Del$. It follows that $g\in G_\a$, a contradiction.

{\bf Case 2.}  Suppose that $K$ is insoluble. In this case, $G_\a$ is intransitive on $\Ome$, and
$K$ has a   normal subgroup $L$ isomorphic to $\A_r$, where $r\in \{5,6,7,8\}$.
Choose a $G_\a$-orbit $\Del'$ such that $L$ is faithful on $\Del'$. Then $m':=|\Del'|\geqslant r$, and $19\geqslant n\geqslant m+m'\geqslant m+r$.

Note that $2^4|H|_2\leqslant |G_\a^{\Del}|_2\leqslant 2^5|H|_2$, and
$|G_\a^{\Del}|_2=|\S_m|_2$ or $|\A_m|_2$.
If $r=8$ then $m\geqslant 12$, and so $n\ge m+r\geqslant 20$, a contradiction.
Suppose  $r=7$. Then $m\geqslant 8$ and $n\geqslant 15$, and so $|G|_2\geqslant 2^{10}$.
It follows that $|G|_2=2^{10}$ and $m=8$; however, in this case, $G_\a^{\Del}\cong 2^4{:}\A_7$, which
can not be contained in a group isomorphic to $\S_8$.
For $r=6$ and $H\cong \A_6$, we get a similar contradiction as above.
Suppose that
$r=6$ and $H\cong \S_6$. Then $2^8\leqslant |G_\a^{\Del}|_2\leqslant 2^9$,
and thus $10\leqslant m \leqslant 13$, yielding $n\geqslant 16$. This leads to
$|G_\a|_2\geqslant 2^{14}$, which is impossible.

By the above argument, we have $r=5$ and $|G_\a|_2=2^{8}$, $2^9$ or $2^{10}$, and then $n\leqslant 15$.
On the other hand,  $2^6\leqslant |G_\a^{\Del}|_2\leqslant 2^8$, we have $m\leqslant 11$, yielding $m=10$ and $n=15$.
It follows that $G=\A_{15}$ and $G_\a=(\Alt(\Del')\times 2^4{:}\S_5)\l\s\t\r$, where
$\s$ is a transposition in $\Sym(\Del')$ and $\t$ is  a product of five disjoint  transpositions in $\Sym(\Del')$.
Then both $G_\a$ and $G_{\a\b}$ have two orbits $\Del'$ and $\Del$ on $\Ome$.
Thus there is no element $g\in \N_G(G_{\a\b})$ such that $\l G_\a,g\r$ is transitive on $\Ome$, a contradiction.
\qed

\subsection{Proof of Theorem~\ref{zthm1}}
Let $G$ be an almost simple group with socle $\A_n$, and let $\Ga$ be $(G,2)$-arc-transitive.

The sufficiency is obvious since the complete graphs $\K_n$ and the odd graphs are clearly 2-arc-transitive under the action of $\A_n$.

The necessity has been established in several lemmas, explained below.
By Lemma~\ref{stab-1}, the vertex stabilizer $G_\a$ is either soluble or  divided into two parts~(1)-(2), according to $G_\a^{\Ga(\a)}$ being almost simple or affine.
For the  case where $G_\a^{\Ga(\a)}$ is almost simple, Lemmas~\ref{main-case}-\ref{PSL(2,9)} show that $\Ga$ is a complete graph or an odd graph.
For the affine case, Lemmas~\ref{soluble-case}-\ref{stab-2} verify the theorem.
\qed


\vskip 30pt

\begin{thebibliography}{99}

\bibitem{Baddeley}
R. W. Baddeley, Two-arc-transitive graphs and twisted wreath products,
{\em J. Algebra Combin}. {\bf 2}(1993), 215-237.



\bibitem{Cameron} P.J. Cameron, {\em  Permutation Groups}, Cambridge University Press,  Cambridge, 1999.

\bibitem{Atlas} J.H. Conway, R.T. Curtis, S.P. Noton, R.A. Parker and
R.A. Wilson, {\em Atlas of Finite Groups}, Clarendon Press, Oxford,
1985. (http://brauer.maths.qmul.ac.uk/Atlas/v3/).

\bibitem{FP2}
X.G. Fang and  C.E. Praeger,
Finite two-arc-transitive graphs admitting a Ree simple group,
{\em Comm. Algebra} {\bf 27(8)} (1999), 3755-3769.

\bibitem{FP1} X.G. Fang and C.E. Praeger, Finite two-arc-transitive
graphs admitting a Suzuki simple group, {\em Comm. Algebra} {\bf
27(8)} (1999), 3727-3754.

\bibitem{Godsil-Royle}
C. Godsil and G. Royle,
{\em Algebraic graph theory}, Springer-Verlag, New York, 2001.

\bibitem{HNP}
A. Hassani, L. Nochefranca and C.E. Praeger,
Two-arc-transitive graphs admitting a two-dimensional projective
linear group, {\em J. Group Theory}  {\bf 2} (1999), 335-353.

\bibitem{IP} A.A. Ivanov and C.E. Praeger, On finite affine
$2$-arc transitive graphs, {\em European J. Combin}. {\bf
14(5)} (1993),
 421-444.

\bibitem{Li1} C.H. Li, On finite s-transitive graphs of odd order,
 {\em J. Combin. Theory Ser. B} {\bf 81} (2001), 307-317.

\bibitem{Li3}
C.H. Li, The finite vertex-primitive and vertex-biprimitive $s$-transitive graphs for
$s\geq4$, {\it Tran. Amer. Math. Soc.} {\bf 353} (2001), 3511-3529.

\bibitem{odd-excep}
C. H. Li, J. J. Li and Z. P. Lu,
Two-arc-transitive graphs of odd order -- I, {\it J. Algebraic Combin.} (to appear).

\bibitem{LSS}
C.H. Li, A. Seress and S.J. Song, s-Arc-transitive graphs and normal subgroups,
{\it J. Algebra} {\bf 421} (2015), 331-348.

\bibitem{Li4}
C.H. Li and Hua Zhang, The finite primitive groups with soluble stabilizers, and the edge-primitive s-arc transitive graphs, {\it Proc. Lond. Math. Soc. (3)} {\bf103} (2011), 441-472.

\bibitem{LS}
M.W. Liebeck and J. Saxl,
The primitive permutation groups of odd degree, \emph{J. London
Math. Soc.} {\bf  (2) 31} (1985), 250-264.

\bibitem{Prag-o'Nan}
C.E. Praeger, An O'Nan-Scott theorem for finite
quasiprimitive permutation groups and an application to $2$-arc
transitive graphs, \emph{J. London Math. Soc.} {\bf 47} (1992),
227-239.

\bibitem{P-W}
C.E. Praeger and J. Wang, On primitive representations of finite alternating and symmetric groups with a $2$-transitive subconstituent, {\em J. Algebra} {\bf 180} (1996), 808-833.

\bibitem{Trofimov}
V.I. Trofimov, Vertex stabilizers of locally projective groups of
automorphisms of graphs. A summary, \emph{Groups, combinatorics and
geometry} (2001), 313-334.

\bibitem{Tutte}
W. T. Tutte, A family of cubical graphs,
\emph{Proc. Cambridge Philos. Soc.} {\bf 43} (1947), 459-474.

\bibitem{Weiss}
R. Weiss, $s$-transitive graphs, Algebraic methods in
graph theory, \emph{Colloq. Soc. Janos Bolyai} {\bf  25} (1981),
827-847.

\end{thebibliography}
\end{document}